%% file: main.tex
\documentclass[journal]{IEEEtran}
\usepackage{graphicx} 

\graphicspath{{images/}}

\usepackage{amsmath}
\usepackage{amsthm}   
\usepackage{amssymb}    
\usepackage{bm}  

\usepackage{dsfont}
\usepackage[dvipsnames]{xcolor}
\usepackage{hyperref}
\usepackage{cleveref}

\usepackage{booktabs}
\usepackage{mathtools} 
\usepackage{xcolor}
\usepackage{url} 
\usepackage{tikz}
\usetikzlibrary{arrows.meta, positioning, shapes.geometric, fit, shadows, automata}

\usepackage{multirow}

\usepackage{caption}
\usepackage{subcaption}

\usepackage{tablefootnote}

\usepackage{dblfloatfix}

\newenvironment{newtext}{\color{black}}{}
\newenvironment{revised}{\color{black}}{}

\title{Electricity Market Bidding for Renewable Electrolyzer Plants: An Opportunity Cost Approach}

\author{\IEEEauthorblockN{Andrea Gloppen Johnsen\IEEEauthorrefmark{1},
Lesia Mitridati\IEEEauthorrefmark{1},
Jalal Kazempour\IEEEauthorrefmark{1}, and
Line Roald\IEEEauthorrefmark{2}} \\
\IEEEauthorblockA{\IEEEauthorrefmark{1}Technical University of Denmark, Kgs. Lyngby, Denmark} \\
\IEEEauthorblockA{\IEEEauthorrefmark{2}University of Wisconsin-Madison, Madison, USA}
\vspace{-6mm}
}

\begin{document}

\input{arxiv_statement}

\maketitle

\vspace{-6mm}
\setcounter{page}{1}

\input{0_Abstract}
\input{1_Introduction}
\input{2_A_Preliminaries}

\input{2_Methodology}
\input{3_Results}

\input{4_Conclusion}

\bibliographystyle{ieeetr}
\bibliography{sample}

\appendix
\counterwithin*{equation}{subsection}
\renewcommand\theequation{\thesection\arabic{equation}}
\input{5_appendix}

\end{document}

%% file: arxiv_statement.tex
\onecolumn
\thispagestyle{empty}

This work has been submitted to the IEEE for possible publication. 
Copyright may be transferred without notice, after which this version 
may no longer be accessible.

\clearpage
\twocolumn

\newpage

%% file: 0_Abstract.tex
\begin{abstract}
Hydrogen produced through electrolysis with renewable power is considered key to decarbonize several hard-to-electrify sectors. This work proposes a novel approach to model the active electricity market participation of co-located renewable energy and electrolyzer plants, based on opportunity-cost bidding. While a renewable energy plant typically has zero marginal cost, selling power to the grid carries a potential opportunity-cost of \textit{not producing hydrogen} when it is co-located with a hydrogen electrolyzer. We first consider only the electrolyzer, and derive its revenue of consuming electricity based on the non-convex hydrogen production curve. We then consider the available renewable energy production and form a piecewise linear cost curve representing the opportunity cost of selling (or revenue from consuming) various levels of electricity. This cost curve can be used to model a stand-alone electrolyzer or a co-located hydrogen and renewable energy plant participating in an electricity market. Our case study analyzes the effects of market-bidding electrolyzers on electricity markets and grid operations. We compare two strategies for a co-located electrolyzer-wind plant; one based on the proposed bid curve and one with a more conventional fixed electrolyzer consumption. The results show that electrolyzers that actively participate in the electricity market lower the average cost of electricity and the amount of curtailed renewable energy in the system compared with a fixed consumption case. However, the difference in total system emissions between the two strategies is insignificant. The specific impacts vary based on electrolyzer capacity and hydrogen price, which determines the location of the co-located plant in the electricity market merit order.

\end{abstract}

%% file: 1_Introduction.tex
\section{Introduction}
The deteriorating effects of climate change demand decarbonization efforts in all sectors. Certain sectors, such as chemical manufacturing, industrial heat and heavy transportation, are difficult to directly electrify. 
However, these sectors can be decarbonized through electrolytic hydrogen production powered by renewable energy sources (RES). 

Electrolytic hydrogen is produced through ``water-splitting", where an electric current is applied to an \textit{electrolyzer}, producing hydrogen (and oxygen) gas from water. The emissions of this process are mainly associated with the consumed power, and are close to zero if renewable power is consumed. The hydrogen can be further processed into a range of chemicals such as ammonia or methanol, which can be used, e.g., for the production of fertilizers or renewable fuels. The world already consumes around 95 Mt of hydrogen gas per year, although this hydrogen is mainly produced through reforming of natural gas, with emission intensities of around 10 kg CO$_2$ eq. per kg H$_2$~\cite{IEA_H2}.  
To decarbonize this hydrogen, governments are planning for the integration of large-scale electrolyzers into the power system, as laid out in the U.S. and E.U. strategies for hydrogen production~\cite{PtX_strategy_USA, PtX_strategy_EU}.

Although hydrogen produced through electrolysis enables low-carbon fuels and chemicals, the build-out of electrolyzers will cause significant increases in power demand~\cite{epri2021}. Theoretically, electrolyzers can serve as flexible and controllable loads, making them well-suited for modern power systems characterized by high shares of RES and limited transmission capacity~\cite{PtX_strategy_DK}. 
\begin{revised} 
Demand-side flexibility and price elasticity are widely recognized as key tools for enhancing the efficiency of electricity markets~\cite{Kirschen_demand_side_2003}. They can also strengthen the business case for RES by mitigating low electricity prices and reducing renewable curtailment.\end{revised} Although electrolyzers are theoretically flexible, it remains unclear whether existing technologies and economic incentives are sufficient to encourage their use of this flexibility in ways that support renewable generation. The primary goal of electrolyzers is to maximize profits from hydrogen production, creating a strong incentive to operate at high utilization rates in order to amortize the substantial capital costs of their construction. Further, electrolyzers might operate with 
\begin{revised} 
fixed-price Power Purchase Agreements (PPAs) and utility contracts that reduce their sensitivity to spikes in the wholesale electricity price~\cite{Pombo-Romero2024}.
\end{revised} 
Electrolyzers might therefore act more like large, constant loads than flexible resources, which could increase the cost of generation and emissions from the power system. 
\begin{revised}
To mitigate these potential drawbacks, particularly the risk of increased emissions, the E.U. and the U.S. have introduced incentive-based regulations aimed at linking hydrogen production to renewable energy generation~\cite{EU_RFNBO, US_IRA}.
\end{revised}

\begin{figure}
    \centering
    \includegraphics[width=\linewidth]{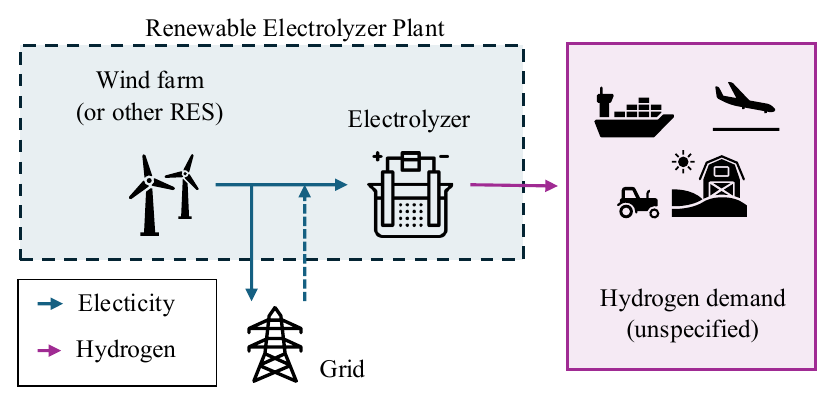}
    \caption{\begin{newtext}\small The renewable electrolyzer plant, REP, consisting of co-located RES and an electrolyzer with a common connection point to the grid. 
    The REP generates electricity from the RES and produces hydrogen via the electrolyzer. The electricity produced can be used internally for hydrogen production or exported to the grid. When RES production is low, the REP may import electricity to maintain hydrogen production. The hydrogen produced is exported to an external demand, which is not specified or further considered here.\end{newtext}}
    \label{fig:REP_illustration}
\end{figure}

This work seeks to explore the potential benefits of \textit{active} electricity market participation by electrolyzers, considering the viewpoint of both the electrolyzers and the grid. We derive cost curves that can model the 
\begin{revised}
price-taking bidding behavior   
\end{revised}of electrolyzers in electricity markets, specifically considering a \emph{Renewable Electrolyzer Plant} (REP) consisting of co-located RES and an electrolyzer with a common grid connection, \begin{newtext}
    as illustrated in Fig.~\ref{fig:REP_illustration}.
\end{newtext} Since the co-location setup reduces the dependency on the transmission grid and allows for direct consumption of renewable power, it is an economically and operationally attractive configuration for large-scale electrolyzer plants. We then examine how the active market participation of large-scale electrolyzers affects both their profitability and the operation of the power system.

\subsection{Literature review}
The market-based scheduling of electrolyzers or REPs is typically modeled as a profit-maximization problem, based on a forecast of day-ahead electricity prices. References \cite{ely_details},~\cite{raheli2023conic},~\cite{ShiYou}, and~\cite{VARELA20219303} 
consider different modeling approaches for scheduling electrolyzers or REPs, with a particular focus on modeling the physical characteristics of the electrolyzer. These references highlight the importance of accurately representing the non-convex \textit{hydrogen production curve}, which describes the relationship between electricity consumption and hydrogen output. They perform optimal scheduling under the assumption of exogenously determined, though typically time-varying, electricity prices, without accounting for the impact of electrolyzer operation on electricity network
constraints and electricity price formation. As electrolyzer capacity expands, their electricity consumption will increasingly contribute to grid congestion and influence electricity prices. 
As a result, assuming that the electricity price is an exogenous signal, as in~\cite{ely_details}–\cite{VARELA20219303}, can lead to infeasible or suboptimal schedules, misleading profit estimates, and negative impacts on power system operations. \begin{newtext}
Reference~\cite{Miljan2022} develops a bi-level optimization model for a hybrid RES-battery-electrolyzer system. This approach endogenizes market prices by removing the assumption of an exogenous price signal, while assuming that the hybrid system participates in the market as a strategic bidder.
\end{newtext}

The impact of electrolyzers on power and energy systems has been explored primarily through capacity expansion models. A technology-agnostic study was conducted in~\cite{jenkins_demand_sinks} employing a generation capacity expansion model to evaluate the role of ``demand-sinks", including electrolyzers, in the grid under various cost assumptions. 
Demand sinks are assumed to be flexible and dispatched according to a fixed utility per consumed energy. References~\cite{US_study_Jenkins,Giovanniello}, and~\cite{elisabeth_zeyen} investigate the emissions associated with grid-connected electrolytic hydrogen production, evaluating the effects of proposed regulations such as temporal matching between electrolyzer consumption and renewable generation, as well as the additionality of matched renewable sources. Reference~\cite{Lueven} explores the impact of various hydrogen investment support policies on the interactions among hydrogen, electricity, and emission markets using a partial equilibrium model. 
However, these studies generally neglect to accurately model the active market participation and price-sensitivity of electrolyzers, focusing instead on simplified assumptions such as fixed loads, which can lead to overlooking the system-wide impact of their operational flexibility.
In addition, these references represent the transmission network either as a zonal system (modeling transmission capacities \textit{between} regions only) or as a copper-plate system (assuming unconstrained transmission capacities). 
Yet, as argued in~\cite{Sofia2024}, ignoring the \textit{deliverability} of renewable energy, by relying on a zonal or copper-plate representation of the transmission network, can significantly distort the analysis of emissions associated with flexible loads such as grid-connected electrolyzers.

\begin{newtext}
    Based on the above literature review, 
    the following research gaps are identified.
    Existing studies that focus on accurate electrolyzer modeling typically rely on either exogenous price assumptions or strategic bidding behavior. However, there is a lack of techno-economic modeling frameworks enabling the accurate representation of REPs into standard market-clearing formulations while still adequately capturing their key physical characteristics. This gap limits the ability to realistically assess the system-wide impact of REPs.
The literature examining the impact of electrolyzers (or REPs) on the power system typically adopts simplified electrolyzer models and overlooks intra-regional transmission constraints. Moreover, while several studies investigate the effects of proposed regulations for electrolytic hydrogen production on system emissions and costs, there is limited analysis comparing the system-level effects of electrolyzers operating under fixed power consumption versus leveraging their inherent flexibility through active market participation.
By deriving the true marginal cost curve of the REP, its scheduling can be directly simulated within electricity market-clearing models, under the assumption of perfect competition, without relying on exogenously determined electricity prices. This cost curve represents the REP as an active market participant and fully captures its price elasticity.
   
Similar marginal cost curves, often represented as price-quantity bidding curves, have been developed in the literature for various assets, or portfolios of assets. For example, the cost curve of thermal generators captures varying marginal production costs~\cite{871739, Conejo2002}, whereas the cost curve of retailers reflects energy procurement costs from the day-ahead market along with balancing costs~\cite{Fleten2005}. For generic portfolios of assets, this model represents opportunity costs arising from coordinated bidding across sequential markets under uncertainty~\cite{Boomsma2014}. One example is the price-quantity bid curve of price-taking wind power producers under uncertainty, accounting for asymmetric imbalance costs~\cite{5373832}. Our proposal follows a similar direction but is novel in its focus on price-taking REPs, whose bidding model in the form of price-quantity curves has not yet been addressed in the literature. 
Unlike previous works that derive bid curves for generic portfolios of assets, our approach explicitly targets the unique techno-economic characteristics of REPs. The proposed bidding curve simultaneously captures the joint physical constraints and flexibility of the co-located renewable and electrolyzer assets, while modeling the cross-sectoral coordination between the power and hydrogen systems, an aspect largely overlooked in existing studies.
%
\end{newtext}



\subsection{Novelty and contribution}
To address the identified gaps in the literature, we make two main contributions.
\begin{revised}
Our first contribution is a practical modeling framework that derives and integrates the bidding models of REPs into market-clearing formulations, providing the basis to derive actionable insights into their system-wide impacts on emissions, costs, and congestion, under perfect competition.  
\end{revised}To do so, we  \textit{analytically} derive the marginal cost curve of REPs or stand-alone electrolyzers, capturing the non-constant efficiencies of the hydrogen production curve and the time-varying availability of renewable generation. The cost curve is then converted to a bid curve that can be used for \begin{revised}
    active participation in, for example, day-ahead electricity markets,
\end{revised}representing either the \emph{opportunity cost} of selling renewable generation to the grid rather than consuming it for hydrogen production, or the \begin{revised}
    \emph{revenue} from hydrogen production derived by importing electricity
\end{revised}when local renewable generation is insufficient. \begin{revised}
Our second contribution builds on this framework by demonstrating its application: using the derived bid curves to derive the market-based scheduling of REPs, we offer a structured method to evaluate their system-wide impacts.
\end{revised}Specifically, we compare a market-bidding REP, which bids based on the derived cost curve, with a fixed-consumption REP, where the electrolyzer operates at a constant load. Our analysis examines how these two distinct strategies affect the total generation cost of the system, renewable power curtailment, and system emissions, highlighting the influence of the REP’s marginal cost and its position within the overall system merit order. Finally, we investigate how transmission constraints impact the perceived effects of connecting an electrolyzer to the power system.

%% file: 2_A_Preliminaries.tex
\begin{newtext}
\section{REP and Market Modeling}
\label{sec:assumptions_market}



In this work, we consider an REP that interacts with a short-term electricity market, as illustrated in  Fig.~\ref{fig:REP_illustration}. 
The REP is connected to the grid through a common connection point shared by both the wind farm and the electrolyzer. From the grid's perspective, the REP appears as a renewable generator with variable power availability. 
However, the presence of the electrolyzer introduces a non-zero opportunity cost of selling power to the grid, as the power could instead be used to produce hydrogen. Further, the REP may become a net importer of electricity when wind production falls short of the electrolyzer’s capacity. Thus, the REP is also comparable to technologies like pumped storage hydropower plants.

The REP trades two products, electric power and hydrogen gas. The electric power is traded through an electricity market, e.g., a day-ahead market. The REP participates in this market through submitting a bid curve, which we later derive. We assume that the hydrogen is traded through a bilateral contract at a fixed price, but we do not specify the kind of contract apart from the fixed hydrogen price. In the following, we outline the main assumptions underlying the derivation of the REP's bid curve and subsequent analysis, and briefly discuss their implications.

%


\subsubsection*{Electricity Market Assumptions}
We derive the REP’s cost curve based on three main assumptions about the electricity market. First, the electricity market is centrally cleared with uniform pricing, which may be either nodal or zonal. Second, the market operates under perfect competition. Third, the electricity market requires submitted bid curves to be piecewise linear and convex (i.e. the price must increase with increasing production). While other bid formats, such as block orders used in E.U. markets, are possible, this assumption reflects the primary bidding format employed in several U.S. and E.U. markets \cite{NEMOCommittee2019, PJM_manual15}.


The first two assumptions imply that the optimal bid curve for any participant, including the REP, is to bid their true marginal cost \cite{Conejo_marginal_pricing}. Our bid curve derivation shows that the REP’s true cost curve must reflect the opportunity cost of selling electricity to the market instead of using it internally for hydrogen production. If the electricity market were cleared using pay-as-bid pricing, for example, the suggested bid curve would not be suitable, as it would not enable the REP to fully recover its investment costs. 


The assumption of perfect competition often does not hold in practice. The REP, which benefits from both low-cost renewable power production and a large controllable load, might possess market power. However, exploiting market power is prohibited in all E.U. electricity markets under REMIT \cite{eu_remit_2011}. We argue that, similar to other generation technologies, the REP's true marginal cost can be audited by regulators or system operators as part of market power mitigation efforts, which supports the validity of the assumption of \textit{marginal cost bidding}.



\subsubsection*{Modeling of Renewable Energy Variability and Hydrogen Demand}
Our derivation of the REP’s cost curve does not consider a specified hydrogen demand (such as a minimum demand over a defined period), since incorporating this would require modeling multiple future time steps and forecasting electricity market-clearing prices. Furthermore, we assume that co-locating the renewable energy source and the electrolyzer within the REP provides sufficient hedging against renewable production variability. As a result, the REP is treated as a fully merchant entity without relying on additional PPAs. However, in our case study, we compare this fully merchant REP, modeled using the derived cost curve, with a REP operating under fixed electrolyzer consumption.



The derived bid curve and case study are both deterministic and do not explicitly model the stochastic nature of the REP or other RES in the system. However, we do account for the variability of RES both in the derivation of the REP bid curve and within the case study.

\subsubsection*{Considered Time Horizon}
The short-term temporal horizon adopted in this study excludes considerations related to optimal investment decisions, such as the size and location of the REP and its long-term profitability. Moreover, potential long-term implications, such as the REP’s influence on future investments in generation and transmission capacity, are not captured within the current framework. In addition, the analysis is limited to a single trading floor and does not consider the balancing market, re-dispatches, other forward markets, or reserve markets.

\subsubsection*{Modeling of Hydrogen Price}
Hydrogen is typically traded through bilateral agreements, such as hydrogen purchase agreements, where the parties negotiate the price \cite{EU_H2_landscape_24}. Therefore, in this work, we assume that the hydrogen price is fixed and exogenously determined. The hydrogen prices used are based on reports of hydrogen production costs, assuming that electrolytic hydrogen competes with conventionally produced hydrogen. The cost of conventionally produced hydrogen via reforming is naturally linked to the price of natural gas. For reference, the International Energy Agency (IEA) reported hydrogen production costs ranging from \$1 to \$4 per kg in 2021, with the upper bound rising above \$8 per kg in 2022 due to increases in natural gas prices \cite{IEA_H2_2024}. We select a lower range value of \$1.5 per kg as the base case, aligning with cost data in the later presented case study from 2019. Additionally, we investigate a scenario where the REP receives a ``green premium'' in the form of a subsidy or production tax credit per kg of hydrogen, based on the highest proposed level of around \$4.5 per kg \cite{EU_RFNBO, US_IRA}. These two prices, i.e., the market rate price of \$1.5 per kg and the subsidized price of \$6 per kg, represent lower and upper bounds of competitive hydrogen prices, excluding impacts such as carbon taxes. While time-varying hydrogen prices would directly influence the REP’s cost curve, we keep the cost of conventional gas generators constant in the case study. Since the market rate hydrogen price correlates with the gas price, a constant hydrogen price is deemed sufficient. Nonetheless, the proposed REP's cost curve framework allows for time-varying hydrogen prices at the same temporal resolution as the electricity market considered.

\end{newtext}

%% file: 2_Methodology.tex
\begin{figure*}[t!]
    \begin{subfigure}[t]{0.24\textwidth}
                \caption{\small  Hydrogen production}
                \label{subfig:REP_CC_A}
        \includegraphics[width = \textwidth]{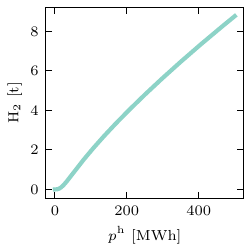}
    \end{subfigure}%
    \hfill
    \begin{subfigure}[t]{0.24\textwidth}
                \caption{\small   Hydrogen revenue}
                                \label{subfig:REP_CC_B}
        \includegraphics[width = \textwidth]{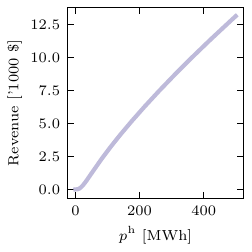}
    \end{subfigure}
        \hfill
    \begin{subfigure}[t]{0.24\textwidth}
                \caption{\small  REP's cost curve}
        \includegraphics[width = \textwidth]{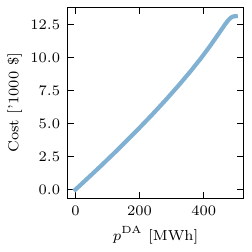}
                        \label{subfig:REP_CC_D}
    \end{subfigure}
    \begin{subfigure}[t]{0.24\textwidth}
                \caption{\small REP's marginal cost curve}
        \includegraphics[width = \textwidth]{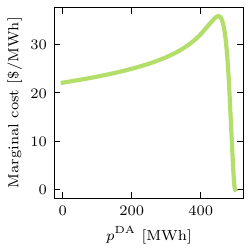}
                        \label{subfig:REP_CC_E}
    \end{subfigure}
\caption{\small Derivation of the REP (marginal) cost curve from the empirical hydrogen production curve. Illustrated for $\bar{P}^{\rm h} = P^{\rm RES} = 500$ MW.}
\label{fig:REP-CC}
\end{figure*}

\section{Methodology}
\begin{revised}
This section begins by outlining the electrolyzer modeling approach used in this work, leading to an empirical hydrogen production curve. Based on this, it derives the \end{revised}REP's continuous, non-convex, (marginal) opportunity cost curve. This cost curve is then approximated as a convex, piecewise linear function suitable for electricity market bidding.  
Finally, the section formulates the electricity market-clearing mechanism as a DC Optimal Power Flow (DC OPF) problem.

\begin{newtext}

\subsection{Modeling of electrolyzers}
\label{sec:assumptions_ely}
This work considers low-temperature electrolyzers, such as proton exchange membrane or alkaline electrolyzers, which are both mature technologies. The electrolyzer produces hydrogen at various efficiencies depending on its partial load set-point, with a peak efficiency typically around 20 to 30 percent of nominal capacity for low-temperature electrolyzers. We base our electrolyzer model on the hydrogen production curve developed in \cite{ely_details}, which does not have a closed form expression. Here, we first express the hydrogen production as a function of power consumption using an empirical model based on the works of \cite{Ulleberg} and \cite{Snchez2018}, resulting in a nonlinear curve. We then approximate it as a piecewise linear, concave function, which is shown to be accurate even for a small number of pieces. The approximation is exact at the linearization points; otherwise it mainly underestimates the original non-concave curve. The largest approximation error is found at the lowest electrolyzer load, the non-concave part of the curve. However, errors at the lower part of the curve have less impact on the scheduling results, as the hydrogen production is naturally the lowest at low loads. Therefore, the approximation error appears to be negligible.


We assume that low-temperature electrolyzers have ramp rates that allow them to freely change their set-point within the time duration of the considered electricity market, which in this case is one hour. Although alkaline electrolyzers are generally slower than proton exchange membrane electrolyzers, ramp rates reaching full load within a minute have been reported for both technologies \cite{Wang2025, Eichman2014}. However, in the case of a full shut-down, alkaline electrolyzers in particular may have start-up times that exceed one hour, which our model does not capture. To avoid full shut-downs, the electrolyzer might instead be operated in a stand-by mode, where it consumes a small amount of electricity (around 2.5\% of nominal capacity), from which fast start up is possible \cite{ely_details}. Our derivation neglects the cost this stand-by consumption would incur. Furthermore, we assume that the electrolyzer plant is operated as a modular system, which enables us to neglect the lower load limit of individual electrolyzer modules \cite{lesniak2024advancedschedulingelectrolyzermodules}.
\end{newtext}

\subsection{Derivation of the REP's  opportunity-cost curve}
To model the REP as an active electricity market participant, we first derive its  cost curve. Typically, renewable generators have zero marginal cost, as there is no variable cost associated with production. However, in a REP that includes both renewable generation and an electrolyzer, the situation is different. 
While producing hydrogen, the REP derives a \textit{revenue} from consuming electricity. Selling this electricity to the grid instead would thus be associated with an \textit{opportunity cost}, leaving the REP with a non-zero cost of exporting electricity.
In the following, we derive a cost curve for the REP based on the opportunity cost associated with lost hydrogen production. The cost curve can be used to model the REP's bids for electricity market-clearing models. From the cost curve, we further derive the REP's \textit{marginal} cost curve. The steps of the derivation are described in detail below and illustrated in Fig.~\ref{fig:REP-CC}.

\subsubsection{Hydrogen production curve}
\label{sec:REP_MC}

The starting point of the derivation is the  \emph{hydrogen production curve} of the electrolyzer,
\begin{align}
    {\rm h}(p^{\rm h}).
\end{align}

This curve expresses the amount of hydrogen ${\rm h}$ produced per hour at a given power consumption $p^{\rm h}$, and does not have a closed form expression. Here, we leverage the empirical hydrogen production curve for an alkaline electrolyzer from \cite{ely_details}, which is shown in Fig.~\ref{subfig:REP_CC_A}. Due to low efficiencies at low and high loads, the hydrogen production curve is non-convex.

\subsubsection{Hydrogen revenue}
Assuming that the produced hydrogen can be sold at a constant price $\lambda^{\rm{h}}$, we derive the \textit{ hydrogen revenue} for a given power consumption ${\rm r}_{\rm h}(p^{\rm h})$ as the product of the hydrogen production and the price of hydrogen,
\begin{align}
\label{eq:hydrogen_revenue}
   &{\rm r}_{\rm h}(p^{\rm h}) = {\rm h}(p^{\rm h}) \lambda^{\rm{h}} & p^{\rm h} \in [0, \bar{P}^{\rm h} ],
\end{align}
where $\bar{P}^{\rm h}$ is the electrolyzer capacity. The resulting hydrogen revenue curve is plotted in Fig.~\ref{subfig:REP_CC_B}.

For an electrolyzer within a REP where the renewable generation capacity is larger than the electrolyzer capacity, i.e., $P^{\rm RES} > \bar{P}^{\rm h}$, the available renewable power may at times exceed the electrolyzer capacity. In these cases, the excess power cannot be used to further increase hydrogen production. 
The maximum hydrogen revenue is therefore equal to the revenue of maximum hydrogen production, ${\rm r}_{\rm h}(\bar{P}^{\rm h})$. To reflect this case, we formulate a piecewise curve, defining the hydrogen revenue for any available (positive) power \begin{revised} $p$:
\begin{align}
\label{eq:hydrogen_revenues}
   {\rm r}_{\rm h}(p) = 
   \begin{cases}
      {\rm h}(p) \lambda^{\rm{h}}, & p \in [0, \bar{P}^{\rm h} ],\\
        {\rm h}(\bar{P}^{\rm h}) \lambda^{\rm{h}}, & p > \bar{P}^{\rm h}.
   \end{cases}
\end{align}
\end{revised}

The first piece \begin{revised}
    of eq.~(\ref{eq:hydrogen_revenues}) 
\end{revised} follows what is illustrated in Fig.~\ref{subfig:REP_CC_A}, and the second piece is a constant equal to the revenue of maximum hydrogen production.

\subsubsection{Opportunity cost}

For the REP to sell power instead of producing hydrogen, the revenue from selling power must exceed the revenue that the REP would earn from withholding power for hydrogen production. The REP's cost of selling a power quantity $p^{\rm DA}$ is thus an opportunity cost, equal to the lost revenue of \textit{not producing} hydrogen from the sold power. \begin{newtext}
The sold power $p^{\rm DA}$ relates to the power consumed by the electrolyzer $p^{\rm h}$ through the internal power balance of the REP:
\begin{align}
    P^{\rm RES} \geq p^{\rm DA} + p^{\rm h}.
\end{align}
\end{newtext}
The opportunity cost ${\rm c}_{\rm el}(p^{\rm DA})$ of selling a power quantity $p^{\rm DA}$ is equal to the difference between the hydrogen revenue of withholding \textit{all} available renewable power and the hydrogen revenue when selling a power $p^{\rm DA}$:
\begin{align}
\label{eq:op_cost_curve}
    {\rm c}_{\rm el}(p^{\rm DA}) & = {\rm r}_{\rm h}(P^{\rm RES}) - {\rm r}_{\rm h}(P^{\rm RES} - p^{\rm DA}), 
\end{align}
\begin{revised}
where ${\rm r}_{\rm h}(.)$ is the hydrogen revenue function defined in eq.~(\ref{eq:hydrogen_revenues}). 
Here, ${\rm r}_{\rm h}(P^{\rm RES})$ represents the revenue from using all available renewable power for hydrogen production, while ${\rm r}_{\rm h}(P^{\rm RES} - p^{\rm DA})$ represents the revenue when a power quantity $p^{\rm DA}$ is sold to the electricity market, since this amount can no longer be used for hydrogen production. 
%
The range of possible $p^{\rm DA}$ depends on the available renewable power relative to the electrolyzer capacity. The REP can always export power $p^{\rm DA}$ up until its available renewable capacity $P^{\rm RES}$.  However, if its available renewable production $P^{\rm RES}$ is less than or equal to  the electrolyzer capacity $\bar{P}^{\rm h}$, i.e., $P^{\rm RES}\leq\bar{P}^{\rm h}$, the REP is also allowed to bid to \emph{import} power. This is expressed mathematically as:
 \end{revised}
\begin{align}
\label{eq:p_DA_def}
    p^{\rm DA} \in \begin{cases}
        [P^{\rm RES}-\bar{P}^{\rm h},P^{\rm RES}], & P^{\rm RES}\leq \bar{P}^{\rm h},\\
        [0,P^{\rm RES}], & P^{\rm RES}>\bar{P}^{\rm h}.
   \end{cases}
\end{align}

In reality, a negative $p^{\rm DA}$ represents a \textit{willingness-to-pay} to import power, and not a cost.  For simplicity, we still refer to the full curve as the (opportunity) cost curve.  

The cost curve, plotted for the case where $P^{\rm RES}=\bar{P}^{\rm h}$ in Fig.~\ref{subfig:REP_CC_D}, is a \textit{shifted reflection} of the hydrogen revenue curve 
\begin{revised}
    ${\rm r}_{\rm h}$(.) over both axes. The vertical reflection arises from the fact that the hydrogen revenue curve is subtracted from a constant ${\rm r}_{\rm h}(P^{\rm RES})$ in the cost curve. The horizontal reflection comes from the fact that the variable $p^{\rm DA}$ appears negatively in the hydrogen revenue function ${\rm r}_{\rm h}(P^{\rm RES} - p^{\rm DA})$. 
\end{revised}
This implies that as the hydrogen revenue curve is (approximately) concave, the REP's cost curve is (approximately) convex. 

In the cases where $P^{\rm RES}$ exceeds the electrolyzer capacity, ${\rm r}_{\rm h}(P^{\rm RES}-p^{\rm DA})$ is equal to the maximum hydrogen revenue ${\rm r}_{\rm h}(P^{\rm RES})$ for $p^{\rm DA} \in [0, P^{\rm RES} - \bar{P}^{\rm h}]$. For this range of sold power, the opportunity cost function thus evaluates to zero, indicating that there is no opportunity cost associated with selling renewable power that exceeds what the electrolyzer can consume. The cost curve then forms a piecewise (non-linear) function. 
The presented method allows for computing the cost curve of a \textit{stand-alone} electrolyzer, by setting $P^{\rm RES}$ to zero. 

\begin{figure*}[t!]
    \centering
    \includegraphics[width = \textwidth]{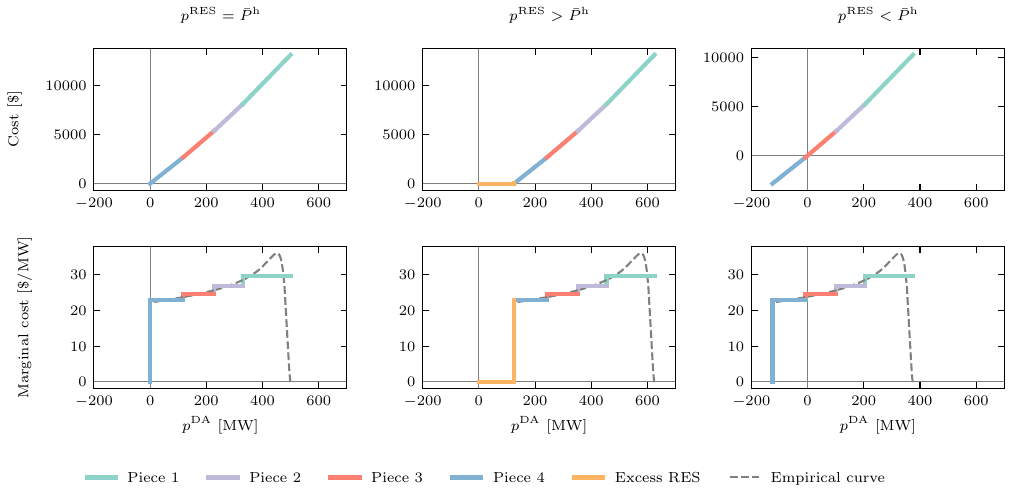}
    \caption{\small Illustrations of the REP's approximated cost curves (first row) and approximated marginal cost curves (second row). We show three cases of renewable power production relative to the electrolyzer capacity, set to 500 MW in this example, with four arbitrarily selected pieces. When there is excess renewable power, meaning production exceeds the electrolyzer capacity, this excess power is offered at zero marginal cost, as it cannot be used for hydrogen production and has no associated opportunity cost. The negative values of the cost curve to the right indicate that the REP bids to purchase the corresponding quantity of power, at a price up to the marginal cost for that quantity.
    }
    \label{fig:MC-cases}
\end{figure*}

\subsubsection{Marginal opportunity cost}
We find the REP's marginal cost of selling power ${\rm mc}_{\rm el}(p^{\rm DA})$ as the change in opportunity cost $d {\rm c}_{\rm el}$ over the change in sold power $d p^{\rm DA}$, i.e., the derivative of the opportunity cost,
\begin{align}
\label{eq:mc}
   {\rm mc}_{\rm el}(p^{\rm DA}) = {\rm c}_{\rm el}'(p^{\rm DA}) = d {\rm c( p^{\rm DA})} / d  p^{\rm DA}.
\end{align}

The marginal cost curve of the REP in the case where $P^{\rm RES} = \bar{P}^{\rm h}$ is illustrated in Fig.~\ref{subfig:REP_CC_E}.

\subsection{Piecewise-linear approximation of the cost curve}
As  illustrated in Figs.~\ref{subfig:REP_CC_D} and \ref{subfig:REP_CC_E}, the REP's (marginal) cost curve is non-convex. However, \begin{revised} as argued in Section~\ref{sec:assumptions_market},
    many
\end{revised}electricity markets require participants to provide piecewise linear, convex cost curves, which corresponds to piecewise constant, non-decreasing marginal costs. 
We therefore derive an approximated cost curve.

Instead of directly approximating the opportunity cost curve in Fig.~\ref{subfig:REP_CC_D}, we begin by approximating the hydrogen production curve in Fig.~\ref{subfig:REP_CC_A}. As will become evident, this allows us to express our approximated opportunity cost curve directly as a function of the hydrogen price $\lambda^{\rm{h}}$, making it easier to adjust in response to changes in the hydrogen price.
We approximate the hydrogen production curve as a concave, piecewise linear function with a decreasing slope, following the approach in \cite{ely_details}, as described in Section~\ref{sec:assumptions_ely}:
\begin{align}
\label{eq:hydrogen_piecewise}
 h(p^{\rm{h}}) = \sum_{i \in \mathcal{I}} \left(A_i p^{\rm{h}} + B_i\right)\mathbf{1}_{\left(p^{\rm{h}} \in [\underline{P}_i, \overline{P}_i ]\right)},
\end{align}
where $A_i$ and $B_i$ represent the slope and constant corresponding to the $i^{\rm th}$ linear piece, and $[\underline{P}_i, \overline{P}_i]$ is the range of electrolyzer power consumption corresponding to piece $i$. The set $\mathcal{I}$ represents the set of pieces \begin{revised}
    $\mathcal{I} = \{1,2,\ldots,I\}$,
\end{revised} and $\mathbf{1}_{(.)}$ is the indicator function which takes on a value of 1 when $p^{\rm{h}} \in [\underline{P}_i, \overline{P}_i ]$ (and zero otherwise). The approximation is chosen such that the slope is decreasing for all pieces, i.e., $A_i>A_{i+1}$. 


From the approximated hydrogen production curve, we derive the approximated opportunity cost for selling $p^{\rm DA}$:
\begin{align}
\label{REPs_OC_approx}
 & {\rm c}_{\rm el}(p^{\rm{DA}})  =  {\rm r}_{\rm h}(P^{\rm RES}) - {\rm r}_{\rm h}(P^{\rm RES}-p^{\rm DA}) \notag  \\ 
  &  = \lambda^{\rm h} \sum_{i \in \mathcal{I}} \left(A_i P^{\rm RES} + B_i\right)\mathbf{1}_{\left(P^{\rm{RES}} \in [\underline{P}_i ,\overline{P}_i ]\right)} 
 \\ 
 & - \lambda^{\rm h} \sum_{i \in \mathcal{I}} \left(A_i (P^{\rm RES}-p^{\rm DA}) + B_i\right)\mathbf{1}_{\left((P^{\rm{RES}}-p^{\rm{DA}}) \in [\underline{P}_i ,\overline{P}_i ]\right)} \notag.
\end{align}

The marginal opportunity cost is the derivative of the opportunity cost, which for the piecewise linear curve is the slope of each piece with respect to $p^{\rm DA}$. The approximated marginal cost therefore becomes an increasing, step-wise function:
\begin{align}
\label{REP_marginal_OC_approx}
 & {\rm mc}_{\rm el}(p^{\rm{DA}}) 
 = \sum_{i \in \mathcal{I}}  \lambda^{\rm h } A_i \mathbf{1}_{\left((P^{\rm{RES}}-p^{\rm{DA}}) \in [\underline{P}_i ,\overline{P}_i ]\right)}.
\end{align}

Since we approximated the hydrogen production curve as a \textit{concave}, piecewise linear function, the resulting cost curve forms a \textit{convex}, piecewise linear function, with increasing slopes. It can therefore be directly implemented into, for example, a DC OPF, without the need for binary variables. 
\begin{newtext}
The slope and intercepts of the resulting piecewise linear cost curve, which form a convex hull, can be expressed through the slope and intercepts of the approximated hydrogen production curve. We achieve a piecewise linear approximation of the REP's cost function ${\rm c_{el}}(p^{\rm DA})$, 
\begin{align}
\label{eq:approx_cost_curve}
{\rm c_{el}}(p^{\rm DA}) = \alpha_i p^{\rm DA} + \beta_i, \quad p^{\rm DA} \in [P^{\rm RES} - \overline{P}_i, \; P^{\rm RES} - \underline{P}_i],
\end{align}
for $i \in \mathcal{I}$ with slopes $\alpha_i$ and intercepts $\beta_i$,
\begin{subequations}
\label{eq:approx_alpha_beta}
\begin{align}
    & \alpha_i = \lambda^{\rm{h}}A_i, \label{eq:approx_alpha} \\
    & \beta_i = \lambda^{\rm{h}}\bigl( (A_j - A_i) P^{\rm RES} + (B_j - B_i) \bigr),
\end{align}
\end{subequations}
where \( j \) is the index such that \( P^{\rm RES} \in [\underline{P}_j, \overline{P}_j] \). Note that the slope in eq.~(\ref{eq:approx_alpha}) corresponds to the marginal cost defined in eq.~(\ref{REP_marginal_OC_approx}). 
In some cases, the RES production may exceed the electrolyzer capacity. To capture this, we define an additional piece on the cost curve. We expand the set of pieces on the approximated cost curve so that  $\mathcal{I} = \{1,2,\ldots, I, I+1\}$ , where $I$ is the number of pieces on the hydrogen production curve. Further, we define that $A_{I+1} =0$, $B_{I+1} =  A_{I} \overline{P}_I + B_{I}$, $\underline{P}_{I+1} =  \bar{P}^{\rm h}$, and $\overline{P}_{I+1} = P^{\rm RES}$. 

The full derivation of the steps from the original formulation of the REP’s cost curve in eq.~(\ref{eq:op_cost_curve}) to the piecewise expressions in eqs.~(\ref{eq:approx_cost_curve}) and~(\ref{eq:approx_alpha_beta}), as well as the definition of the additional piece $I+1$ is provided in Appendix~\ref{app:cost_curve}.


\end{newtext}

Depending on the electrolyzer capacity relative to the available renewable power, three different cases of the (marginal) cost curve arise, as illustrated in Fig.~\ref{fig:MC-cases}:
\begin{itemize}
    \item If the renewable power production is equal to the electrolyzer capacity, $P^{\rm RES} = \bar{P}^{\rm h}$, any sold power is associated with a (non-zero) positive opportunity cost, forming a cost curve that intercepts with the origin. This case is illustrated by the first column of plots in Fig.~\ref{fig:MC-cases}, and corresponds to the approximation of the curves in Fig.~\ref{fig:REP-CC}.

    \item 
    If there is \emph{excess} renewable power production, $P^{\rm RES}>\bar{P}^{\rm h}$, the cost curves are shifted right to higher values of $p^{\rm DA}$.  For the excess renewable generation, \begin{revised}
        the additional piece defined evaluates to zero, as there is no opportunity cost of hydrogen production associated.
    \end{revised}This case is illustrated in the second column of graphs in Fig.~\ref{fig:MC-cases}.
    
    \item If there is \emph{insufficient} renewable power, $P^{\rm RES} < \bar{P}^{\rm h}$, the marginal cost curve is shifted left towards lower values  of $p^{\rm DA}$, and down towards lower costs. 
    The negative bid curve at negative values of $p^{\rm DA}$ is interpreted as the marginal bid for ``negative exports'', which corresponds to imports, i.e., consumption, of electricity. This case is illustrated by the plots in the last column of Fig.~\ref{fig:MC-cases}. 
\end{itemize}

\subsection{Integrating Electrolyzers in the Electricity Market Clearing} 
\label{opportunity_cost_bidding_section}
The market-bidding REP submits its approximated cost curve, derived in the previous section, to the market operator and receives a generation schedule  $p^{\rm DA}$ after the electricity market has cleared. 
The REP then uses the remaining power for hydrogen production (up to the full capacity of the electrolyzer), i.e., $p^{\rm h}= \min \left( P^{\rm RES}-p^{\rm DA}, \bar{P}^{\rm h}\right)$.

To model the electricity market clearing, we solve a standard DC OPF. We assume a system with a set of $\mathcal{N}$ nodes, indexed by $n$. For simplicity of notation, but without loss of generality, we assume that we have one generator $p_n^{\rm G}$ and one inflexible demand $P_n^{\rm D}$ per node. The cost of generator $n$ is modeled using piecewise linear costs curves with a set of $\mathcal{I}_n$ pieces and the auxiliary variable $c_{n}^{\rm G}$ representing the total cost. With this notation, we formulate the DC OPF as follows: 

\begin{subequations}
\allowdisplaybreaks
\label{eq:DC OPF}
{\small  
\begin{align}
\allowdisplaybreaks
  \underset{c_{n}^{\rm G} , p_i^{\rm G}, \theta_{i}}{ \operatorname{min} } & \sum_{n \in \mathcal{N}} c_{n}^{\rm G}  \label{DC_obj} \\
 \text { s.t. }  &  c_n^{\rm G} \geq \alpha_i p_n^{\rm G} + \beta_i && \!\!\forall i \in \mathcal{I}_n, n \in \mathcal{N} \label{PM_gen}  \\
& p_{n}^{\rm G}\!-\!P_n^{\rm D} =\!\!\!\!\! \sum_{m \in \mathcal{K}_n} \!\frac{1}{X_{n m}} \left(\theta_{n}-\theta_{m}\right) && \!\!\forall n \in \mathcal{N} ~~ : \lambda^{\rm DA}_n \label{PM_pb}\\
& \!\!\!\!\!-\!\bar{P}_{nm} \leq  \! \frac{1}{X_{nm}}\left(\theta_{n}-\theta_{m}\right) \leq \bar{P}_{nm} && \!\!\forall m \in \mathcal{K}_n, n \in \mathcal{N}  \label{PM_flowlimits} \\
& \underline{P}_n^{\rm G} \leq p_{n}^{\rm G} \leq \bar{P}_n^{\rm G}  && \!\!\forall n \in \mathcal{N} \label{PM_genlimits}\\
& \theta_{n = 1}=0,  \label{PM_refnode}
\end{align}
}%
\end{subequations}
where the variables $c_n^{\rm G}$, $p_n^{\rm G}$ and $\theta_{n}$ denote the total generation cost, power generation schedule, and voltage angle at node $n$, respectively. The objective function \eqref{DC_obj} minimizes the total cost of generation. Constraint \eqref{PM_gen} implements the piecewise linear, convex cost curves with parameters $\alpha_i$ and $\beta_i$, representing the slope and intercepts for each piece. Note that the inequality will always be binding (i.e., satisfied with equality) as the cost variables $c_n^{\rm G}$ are minimized.  The reactance of each branch is given by ${X_{nm}}$, where $m \in \mathcal{K}_n$ notes all nodes $m$ connected to  $n$. Constraint (\ref{PM_pb}) enforces the nodal power balance, while (\ref{PM_flowlimits}) and (\ref{PM_genlimits}) ensure minimum and maximum limits on branch power flow and generation capacity. Finally, (\ref{PM_refnode}) defines the reference node. The Locational Marginal Prices (LMPs) in the  day-ahead stage are given by the dual variable associated with the power balance constraints (\ref{PM_pb}), denoted as $\lambda^{\rm DA}_n$.
REPs are modeled as generators in this formulation, with the day-ahead schedule $p^{\rm DA}$ represented by a generator variable $p_n^{\rm G}$. Here, we use the slope and intercepts as derived in {\eqref{REPs_OC_approx}} to represent the REP. However, assuming no start-up costs, the approximated marginal cost curve is equivalent to the full cost curve and can also be implemented in a DC-OPF. Further, the marginal cost curve  is important for understanding the REP's position in the merit order curve.

Note that cases where the REP imports power are accommodated by adapting the generation limits to allow for negative generation. We consider no inter-temporal constraints such as  ramping rates or start-up time. 

\begin{newtext}

\input{images/schematic}

Fig.~\ref{fig:schematic} summarizes how the REP is modeled as a market participant in this work.
The REP receives an exogenous hydrogen price, $\lambda^{\rm h}$. This price, together with the (forecasted) internal RES production at time $t$, $P^{\rm RES}_t$, and the electrolyzer's hydrogen production curve ${\rm h}(p^{\rm h})$, is used to derive the REP's piecewise linear bid curve. This bid curve is given by slopes $\alpha_{it}$ and intercepts $\beta_{it}$ for piece $i$ on the curve. The bid curve is submitted to the electricity market, here modeled as a DC OPF (\ref{eq:DC OPF}), which outputs the cleared schedules $p^{\rm G*}_{nt}$ and prices $\lambda^{\rm DA*}_{nt}$ of all generators and nodes in the system.
\end{newtext}

%% file: images/schematic.tex
\definecolor{cb_green}{HTML}{b3de69}
\definecolor{cb_blue}{HTML}{80b1d3}
\definecolor{cb_turquoise}{HTML}{8dd3c7}
\definecolor{cb_red}{HTML}{fb8072}
\definecolor{cb_lilac}{HTML}{bebada}
\definecolor{cb_purple}{HTML}{bc80bd}

\begin{figure}[t]
\centering
\begin{tikzpicture}[
    every node/.style={font=\small},
    source/.style={rectangle, 
    draw=cb_green, fill=cb_green!15,
    thick, minimum width=2cm, minimum height=1.1cm, rounded corners},
    load/.style={rectangle, 
    draw=cb_blue, fill=cb_blue!15,
    thick, minimum width=2cm, minimum height=1.1cm, rounded corners},
    market/.style={rectangle, draw=black, fill=gray!10, thick, minimum width=2cm, minimum height=1.2cm, rounded corners},
    boxgen/.style={rectangle, double copy shadow, draw=black, fill=red!10, minimum width=0.9cm, minimum height=0.9cm, thick, rounded corners},
    boxload/.style={rectangle, double copy shadow={shadow yshift=-.5ex}, draw=black, fill=cyan!10, minimum width=0.9cm, minimum height=0.9cm, thick, rounded corners},
    arrow/.style={-Latex, thick},
    feedback/.style={dashed, -Latex, thick},
    repbox/.style={draw=black, thick, dashed, inner sep=0.3cm, rounded corners},
    node distance=1.6cm and 3cm 
]

\node[source, align=center] (res) at (0,0) {\textbf{RES} \\  $P^{\mathrm{RES}}_{t}$};
\node[load, align=center] (flex) at (0,-1.5) {\textbf{Electrolyzer}\\${\rm h}( p^{\mathrm{h}})$};
\node[market, align = center] (market) at (3.7,-0.75) {Electricity \\ Market \\ (\ref{eq:DC OPF}) };


\node[repbox, fit=(res)(flex), label=above:{\textbf{REP}}, inner sep = 5pt] (rep) {};

\draw[arrow, thick] (rep.east) -- (market.west) node[midway, above] {$\alpha_{it}, \beta_{it}$};
\draw[arrow, thick] (rep.east) -- (market.west) node[midway, below] {(\ref{eq:approx_alpha_beta})};

\draw[-Latex, thick]([yshift=0.7cm]market.north) node[above] {Other generators' cost curves} -- (market.north)           
;

\draw[-Latex, thick]([yshift=-0.7cm]market.south) node[below] {$P^{\rm D}_{nt}$} -- (market.south)           
;
\draw[Latex-, thick] (rep.west) -- ++(-0.7,0) node[left] {$\lambda^{\rm h}$};

\draw[arrow, thick] (market.east) ++(0,0) -- ++(0.7,0) node[right, align=center] {$\lambda^{\rm DA*}_{nt}$\\ \\$p^{\rm G*}_{nt}$};


\end{tikzpicture}
\caption{\small 
\begin{newtext}
Schematic of an REP actively participating in an electricity market (here, a day-ahead market). An exogenous, fixed hydrogen price $\lambda^{\rm h}$, combined with known local renewable power production $P^{\rm RES}_t$ at hour $t$ and hydrogen production curve ${\rm h}( p^{\mathrm{h}})$, defines the REP’s bid curve for that hour. This  curve is characterized by piecewise slopes $\alpha_{it}$ and intercepts $\beta_{it}$, where $i$ indexes the segments. The bid is submitted to the electricity market, together with the cost curves of other generators and the inflexible demand at the time $P^{\rm D}_{nt}$, which returns the clearing prices $\lambda^{\rm DA*}_{nt}$ and production schedules $p^{\rm G*}_{nt}$.
\end{newtext}
}
\label{fig:schematic}
\end{figure}

%% file: 3_Results.tex
\section{Case study}
This section analyzes the impact of connecting an electrolyzer to the grid, when co-located with a renewable power plant, comparing active market participation through the proposed bidding approach to a traditional fixed electrolyzer consumption.

\subsection{Case Study Description}

\begin{figure}[t!]
    \centering
\includegraphics[width = \linewidth]{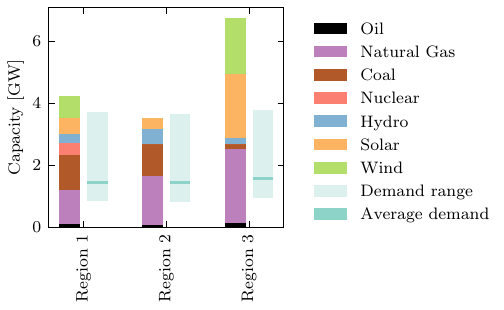}
    \caption{\small Generation capacities and demand range per region of the RTS-GMLC. While region 3 has the most capacity, it also has a slightly higher demand load, approximately 10\% higher than the other two regions across the year. Further, the increase in generation capacity in region 3 is dominated by renewable power sources. For conventional sources, the generation capacity is similar in each region. 
}
    \label{fig:reginal_capacities}
\end{figure}

We define our case study based on the Reliability Test System---Grid Modernization Lab Consortium (RTS-GMLC) \cite{GMLC}, an updated version of the IEEE RTS-96 with high shares of renewable power capacity. The RTS-GMLC system is divided into three regions, and comes with
one year of data for hourly 
load and renewable production, as well as piecewise cost data for all conventional generators. 
The load in the system is inflexible.  We refer to the RTS-GMLC system without any modifications, i.e., without any electrolyzer capacity, as the \emph{base system}. The production capacities and demand range of each region in the system are plotted in Fig.~\ref{fig:reginal_capacities}. Note that the majority of the renewable capacity in the system is located in region 3.

We select a wind farm, generator no. 156 located at node 303 in region 3, to convert to a REP. This wind farm has a total installed capacity of 847 MW and an average production of 237 MW throughout the year. We consider three different electrolyzer capacities, namely 100 MW, 500 MW and 1000 MW, but keep the REP's renewable capacity equal to the pre-existing capacity in all cases. Hence, the only change to the base system is the added electrolyzer capacity at node 303.

We consider two different hydrogen prices, \$1.5 and \$6 per kg. These prices are representative of the current market rate for conventionally produced hydrogen and the current market rate for green hydrogen with a subsidy \cite{PtX_strategy_USA}, respectively. For every hour, we calculate the marginal cost curve of the REP based on the production of the local wind farm, the capacity of the electrolyzer, and the assumed hydrogen price.

\subsection{Analysis of the Merit Order Curve}
\begin{figure*}[!t]
    \centering
\includegraphics[width = \textwidth]{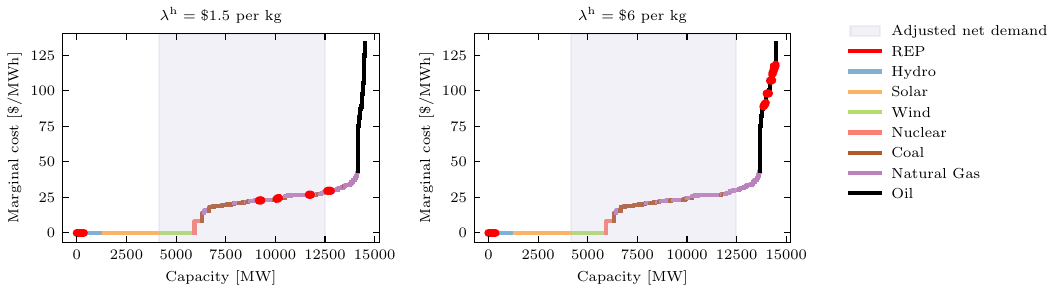}
    \caption{\small Merit order curve of the total installed generation capacities for a low (left) and high (right) hydrogen price. The REP illustrated consists of a 847 MW wind farm and a 500 MW electrolyzer capacity. The adjusted net demand, i.e., the demand subtracted by the renewable production of every time-step, adjusted by adding the total renewable capacity in the system, is indicated by the shaded area. When minimizing cost of generation, only the generation to the left of the adjusted net demand will be dispatched, unless there is congestion in the network. 
    }
    \label{fig:merit_order}
\end{figure*}

To provide an insight into how the REP's marginal cost of electricity generation compares with other generators in the system, we first analyze where the REP falls in the merit order. Fig.~\ref{fig:merit_order} shows the merit order curve of the generation stack of the entire RTS-GMLS system for the two different hydrogen prices, both with an electrolyzer capacity at the REP of 500 MW. The range of \textit{adjusted net-demand} is marked by the purple shaded area. The adjusted net-demand of an hour is the total demand minus the total renewable production, adjusted by adding the total renewable capacity of the system. The range indicates all possible values of adjusted net demand throughout the year. In the absence of congestion, generators that fall in this range would be the marginal generator at certain hours. Generation capacity that falls to the right of the adjusted net demand range will only be dispatched in the event of network congestion.

For both the lower and higher hydrogen prices, the REP (marked in red) has a piece of its marginal cost curve equal to zero, as its renewable capacity of 847 MW exceeds its electrolyzer capacity of 500 MW. However, the presence and size of the zero marginal cost piece will depend on the renewable power availability of every hour. For the lower hydrogen price (left plot), the remaining pieces of the REP cost curve have marginal cost similar to that of natural gas generation and fall partly within the adjusted net demand range. For the higher hydrogen price  (right plot), the remaining pieces of the REP cost curve have higher marginal costs, similar to that of oil-fueled generation. These pieces fall outside the adjusted net demand range, indicating that this power will only be exported when network congestion causes prices to increase.

\subsection{Simulation and Benchmarking Framework}
To analyze the impact of a market-bidding REP on the system and benchmark its effects against a more traditional style of operation, we simulate the electricity market clearing by solving the DC OPF problem \eqref{eq:DC OPF} for each hour of the year. 
To solve the DC OPF, we use the implementation provided in the Julia package \textit{PowerModels.jl} \cite{8442948}. Our source codes are publicly available in \cite{our_git}.

\subsubsection{REP Representations} We compare our proposed marginal cost bidding approach with a more traditional fixed consumption case. Below, we describe each case in detail: 

\textbf{Market-bidding REP:} The REP participates in the market, submitting its opportunity cost curve. A DC OPF is solved for the RTS-GMLC with the cost curve of wind farm 156 altered to represent the presence of the electrolyzer as described in Section~\ref{opportunity_cost_bidding_section}. This case represents how the REP would behave when participating in the electricity market and fully exposed to the LMPs of the system (i.e., fully flexible).

\textbf{Fixed consumption REP:} We compare the market-bidding REP to an REP where the electrolyzer is consuming at a fixed, constant power. To create a fair comparison, we define the fixed consumption power as the average consumption of the market-bidding REP, i.e., the total electricity consumption of the market-bidding REP over the year, divided by the number of hours in the year. The electrolyzer is modeled as an additional load at the REP's node, not modifying the bid curve of wind farm 156. This new load is set to the calculated fixed consumption and is inflexible.
However, for several hours, this level of consumption causes the DC OPF to become infeasible. In those hours, we first model the fixed consumption as a bid with a high marginal cost (of \$10,000 per MWh), thus approximating the highest electrolyzer consumption possible while maintaining a feasible program. We then rerun the DC OPF with the new load set equal to this (lower) consumption.

\subsection{System-focused Results}
\label{sec:sys_results}
\begin{figure*}[!t]
        \centering
\includegraphics[width =\textwidth]{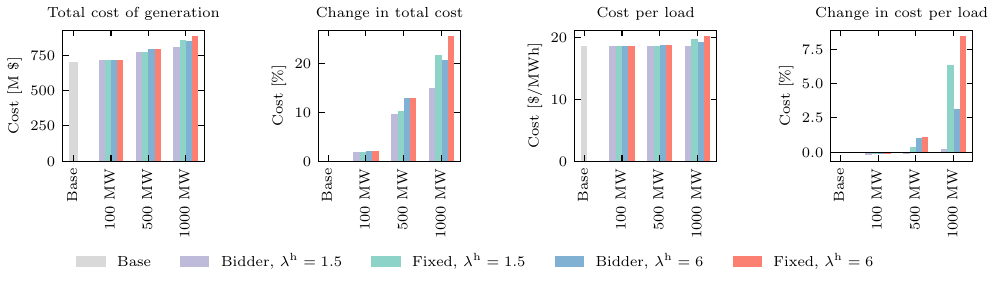}
    \caption{\small Total cost of generation and cost of generation per load for a fixed consumption and market-bidding electrolyzer under the hydrogen prices of $\$$1.5 and $\$$6  per kg.
}
    \label{fig:sys_costs}
\end{figure*}
We first investigate the impact of the REP's electrolyzer from a system perspective, considering the cost of generation, renewable power curtailment, and emissions. We combine the three electrolyzer capacities (100, 500 and 1000 MW), two hydrogen prices (\$1.5 and \$6 per kg) and two REP representations previously described, resulting in 12 cases.

\subsubsection{Cost of generation} 
\label{sec:system_cost}
We examine the average cost of generation throughout the year for each case. Due to the REP's ability to act as a consumer, the objective value of the DC OPF consists not only of the cost of generation, but also of the utility of consumption from the REP. To evaluate the changes in cost of generation due to the presence of the electrolyzer, we compute the cost of generation of each hour as:
\begin{align}
\label{}
    \text{Cost of generation} =  \sum_{g \in \mathcal{G \backslash \text{REP} }} c_{g}^{\rm G},
\end{align}
where $\mathcal{G \backslash \text{REP}}$ is the set of all generators apart from the REP. The cost of generation is therefore \textit{not} equal to the objective value (\ref{DC_obj}). We sum across all hours to achieve the \textit{total} cost of generation, and divide by the total generation in the system to obtain the cost of generation \textit{per load}. The results are shown in Fig.~\ref{fig:sys_costs}, which also reports changes in cost relative to the base system. 

We observe that the total cost of generation (two left plots) increases as the electrolyzer capacity increases. This is expected as the electrolyzer imposes an additional load on the system.
For the smaller electrolyzer capacity of 100 MW, the cost increases by less than 2.5\% for any operational mode or hydrogen price. However, for larger electrolyzer capacities, the total cost of generation increases by 10-20\%.  In general, the market-bidding REP yields a smaller increase in the total cost of generation compared to the corresponding fixed consumption case. This is expected as the market-bidding REP will reduce its consumption if the LMP is sufficiently high. 


Next, we consider the cost of generation per load, shown in the two plots to the right in Fig.~\ref{fig:sys_costs}. We first observe that adding a small electrolyzer with 100 MW capacity decreases the cost of generation per load for both hydrogen prices and REP representations. This is because part of added REP consumption is met by otherwise curtailed renewable power. For electrolyzer capacities of 500 and 1000 MW, the average cost generally increases. However, the magnitude of the increase depends on the case. Cases with fixed consumption and/or a higher hydrogen price lead to significant increases in average cost, reaching up to 7.7\% (although a market-bidding REP causes smaller increases compared to a fixed consumption REP). In contrast, adding a market-bidding REP with a lower hydrogen price yields a negligible increase in generation cost per load.

From these results, we conclude that to keep the cost of generation low, the REP must be sensitive to the LMP, i.e., operated as a market bidder. Further, higher hydrogen prices reduce the benefits of market-bidding REP operation. Following the merit order curve in Fig.~\ref{fig:merit_order}, the market-bidding REP with a high marginal cost is only dispatched if there is network congestion, and will therefore impact the total cost of generation in a way similar to  that of a fixed consumption REP. 
We observe this effect in the cases with a hydrogen price of \$6 per kg and large electrolyzer capacities, and note that subsidies for hydrogen production, through inflating the hydrogen price, might result in less flexible operation of the REP and thus increase the cost of generation of its connected system.

\begin{figure}[t]
        \centering
\includegraphics[]{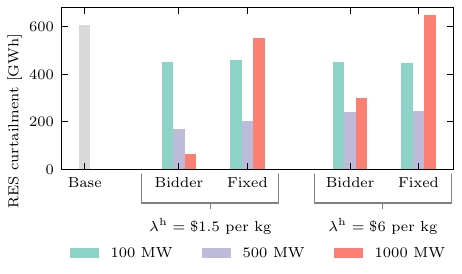}
    \caption{\small Total curtailment of renewable power generation.}
    \label{fig:curtailment_system}
\end{figure}

\subsubsection{Renewable power curtailment}
\label{sec:res_curtailment}
We next investigate how the presence of an electrolyzer in the system affects the curtailment of renewable power, as it is expected that electrolyzers can use surplus renewable power and thus reduce the need for curtailment. 
Fig.~\ref{fig:curtailment_system} shows the total renewable power curtailment over the year in the entire system, for the different electrolyzer capacities, REP representations, and hydrogen prices explored previously. 
\begin{revised}
    In absence of congestion, increasing system load should yield decreasing renewable curtailment due to the renewable generators' position on the merit order curve of Fig.~\ref{fig:merit_order}.
\end{revised} This trend is observed in the market-bidding cases with a low hydrogen price. However, in the other cases, increasing the electrolyzer capacity from 500 to 1000 MW actually increases the curtailment. 
This increase in curtailment is a result of congestion in the network. When the REP bids to import power at a sufficiently high price (or in the fixed consumption cases), congestion forces the DC OPF to dispatch expensive, non-renewable generation and curtails renewable generation to allow the REP to import more power.

We draw a conclusion similar to that in Section~\ref{sec:system_cost}, that high hydrogen prices reduce the flexibility of electrolyzers, potentially leading to negative impacts on RES curtailment.

\subsubsection{System emissions}
Finally, we investigate the emissions associated with the generation of power in the system. 
We compute the total emissions and emissions per load, i.e., total emissions divided by total system load across all hours, as well as the increase from the base system for the 12 described cases. The results, reported in Table~\ref{tab:sys_emissions}, are given in CO$_2$ equivalents. We assume constant emission factors of 0.9606, 0.6042, and 0.7434 metric tons of CO$_2$ per MWh for coal, natural gas, and oil-fueled power plants, respectively, and zero emissions for renewable and nuclear generation. 

Our results show that for a given electrolyzer capacity, both the total emissions and emissions per load remain similar irrespective of both how the REP is operated and the hydrogen price. \begin{revised}
    Hence, there is no clear difference in the additional emissions from the system depending on the operational strategy of the REP, but rather depending on its overall consumption.
\end{revised} This result is understood by examining the merit order curve of the system in Fig.~\ref{fig:merit_order}. The cheapest generation sources, i.e., renewable and nuclear capacities, are indeed also the least emission intense. However, the (generally) cheapest conventional generation capacities, the coal-fueled plants, have higher emission factors than the more expensive natural gas and oil-fueled plants. \begin{revised}
    Therefore, the market-bidding REP, while following price-signals from the electricity market, does not guarantee lower system emissions compared to the fixed REP.
\end{revised}

The emissions generally increase with the installed electrolyzer capacity. For the 100 MW electrolyzer, this increase is negligible, as about 20\% of the additional load the electrolyzer imposes on the system is met with otherwise curtailed renewable power. 
For larger electrolyzer capacities, conventional generation meets a larger share of the additional system load, resulting in increased system emissions.

We conclude that market bidding is only an efficient measure to avoid increases in system emissions under strong correlation between the marginal cost and emissions of generation. 

\begin{table}[!t]
\centering
\caption{\small System emissions results over the year.}
\label{tab:sys_emissions}
\begin{tabular}{lllrrrr} \toprule
   & & & \multicolumn{2}{c}{\textbf{Total   emissions}} & \multicolumn{2}{c}{\textbf{Emissions per   load}} \\
                         \multicolumn{3}{c}{\textbf{Case}}         & {[}M t{]} & $\Delta$ & {[}kg/MWh{]} & $\Delta$ \\ \toprule
\multicolumn{3}{c}{Base system}                      & 13.88 &  \multicolumn{1}{c}{-}   & 368.50 &   \multicolumn{1}{c}{-}       \\ \midrule
\multirow{4}{*}{\rotatebox[origin=c]{90}{100 MW}}  & \multirow{2}{*}{$\$$1.5/kg} & bidder & 14.38 & 3.61\%  & 373.83 & 1.45\%  \\
                         &                             & fixed  & 14.38 & 3.62\%  & 373.77 & 1.43\%  \\
                         & \multirow{2}{*}{$\$$6/kg}   & bidder & 14.43 & 3.97\%  & 374.39 & 1.60\%  \\
                         &                             & fixed  & 14.43 & 3.97\%  & 374.39 & 1.60\%  \\  \midrule
\multirow{4}{*}{\rotatebox[origin=c]{90}{500 MW}}  & \multirow{2}{*}{$\$$1.5/kg} & bidder & 16.39 & 18.13\% & 396.37 & 7.56\%  \\
                         &                             & fixed  & 16.39 & 18.10\% & 396.15 & 7.50\%  \\
                         & \multirow{2}{*}{$\$$6/kg}   & bidder & 16.88 & 21.62\% & 401.42 & 8.93\%  \\
                         &                             & fixed  & 16.88 & 21.64\% & 401.50 & 8.96\%  \\ \midrule
\multirow{4}{*}{\rotatebox[origin=c]{90}{1000 MW}} & \multirow{2}{*}{$\$$1.5/kg} & bidder & 17.70 & 27.54\% & 410.67 & 11.44\% \\
                         &                             & fixed  & 17.69 & 27.46\% & 409.53 & 11.13\% \\
                         & \multirow{2}{*}{$\$$6/kg}   & bidder & 18.08 & 30.29\% & 414.96 & 12.61\% \\
                         &                             & fixed  & 18.32 & 32.00\% & 415.68 & 12.80\% \\ \bottomrule
\end{tabular}
\end{table}

\subsection{REP-focused results}
\label{sec:REP_results}

\input{images/REP_results_table}

We now analyze our results from the perspective of the REP. For each of the 12 cases described previously, we calculate the average load of the electrolyzer, the REP's total profit and the average LMP at the REP's node. The total profit is calculated considering the electricity cost and revenue (assuming that the REP is remunerated for exports or pays for imports based on the hourly nodal LMP) as well as the revenue from hydrogen production (using the hydrogen price for the specific case). 
The results for the 12 REP cases are provided in Table~\ref{tab:REP_results}.

\subsubsection{Electrolyzer utilization}
\label{sec:ely_utilization}
We first discuss the average load of the electrolyzer in the various cases, shown in the first column of Table~\ref{tab:REP_results}. For the cases of the higher hydrogen price and installed electrolyzer capacity of 100 MW and 500 MW, the average load is 100 MW and $>$499 MW, indicating that the overall utilization is (near) 100\%. For the same electrolyzer capacities, but at the lower hydrogen price, the average load is reduced to an overall utilization of 91\% and 84\%. 
These results align well with our analysis of the merit order curve in Fig.~\ref{fig:merit_order}. For a higher hydrogen price, the REP lies high in the merit order curve, and it would seldom be beneficial to reduce consumption to sell electricity. For a lower hydrogen price, the marginal opportunity cost is lower, and selling electricity instead of producing hydrogen is cost-optimal in more hours. Similarly, the lower marginal revenue associated with a lower hydrogen price would cause the REP to import less electricity. 

For an electrolyzer with 1000 MW capacity, we observe that the average load is significantly below capacity for the market-bidding REPs, even for the case with a high hydrogen price. This indicates that a combination of congestion and high prices at the REP node reduces the viability of importing sufficient power to run the electrolyzer at full capacity. 

Further, while the fixed consumption REP should consume the average electrolyzer load of the corresponding market-bidding case, we observe that the actual average consumption is significantly lower than for the market-bidding REP. As previously explained, the fixed consumption level leads to an infeasible program in certain hours, and for those hours we reduce the consumption level of the fixed electrolyzer. The lower average load therefore indicates that a fixed consumption at this level is frequently infeasible.

\subsubsection{Profit}
We next discuss the REP's operational profit under the different cases, shown in the second to last column of Table~\ref{tab:REP_results}.
Recall that the REP profit is calculated assuming that the REP imports and exports power at a price equal to the LMP for all cases. Under this assumption, we hypothesize that the market-bidding REP will achieve a profit equal to or higher than that of the fixed consumption REP, since the market-bidding REP is actively adapted to the LMP. From our results, we observe that this hypothesis holds true. For the smaller electrolyzer capacity and high hydrogen price, the profit is similar for market-bidding and fixed consumption operation, as it is cost-effective to run the electrolyzer at full capacity in both cases. However, if the hydrogen price is lower and/or the electrolyzer capacity is sufficiently large (e.g., 500 or 1000 MW) to cause congestion or infeasibility, we observe that there can be significant differences in profit between the market-bidding and fixed consumption cases.
In the case of a 500 MW electrolyzer capacity and a lower hydrogen price, the market-bidding case yields a 26\% profit increase compared to the fixed consumption case. 
In the case of a 1000 MW electrolyzer capacity, the profits achieved through active market bidding far exceed those of a fixed consumption, which both exhibit a loss (i.e., a negative profit).

\begin{newtext}
While the market-bidding REP guarantees non-negative operational profits, it must still cover its investment costs over the long term. Thus, long-term profitability depends on operational profits being sufficient to recover the costs of constructing the wind farm and electrolyzer over the plant’s lifetime. Since this study focuses on the short term and does not consider factors such as optimal sizing and placement of REP components, we do not analyze long-term profitability. However, our results indicate a tendency for operational profits to decrease as electrolyzer capacity increases, which may serve as an important signal for investment decisions.
\end{newtext}

When we compare profits for market-bidding REPs with 500 and 1000 MW electrolyzer capacity, we observe that a larger capacity increases profits in the case with a \$1.5 per kg hydrogen price but decreases profits in the case with a \$6 per kg hydrogen price. The average LMP at the REP's node, given by the last column of Table~\ref{tab:REP_results}, shows that there is a significant increase in LMP when the electrolyzer capacity increases from 500 to 1000 MW in certain cases. We examine the case of a market-bidding REP under the higher hydrogen price. In this case, the LMP increases from \$31.77 to \$81.54 per MWh. Recalling the merit order curve of Fig.~\ref{fig:merit_order}, the marginal cost of the REP under a higher hydrogen price lies in the range of \$80-120 per MWh. The REP therefore makes more profit from producing hydrogen than from selling electricity at LMPs up to this range. However, with increasing LMP, the profit associated with hydrogen production decreases. The increase in electrolyzer capacity when operating under the higher hydrogen price increases the LMP to a point where the REPs profit decreases.

\begin{newtext}
    Typically, a generator’s marginal cost curve reflects the electricity price at which it breaks even for a given output. Hence, if the cleared price is above its marginal cost, it makes a profit. Essential to the REP is that when wind production is available, the marginal \textit{opportunity} cost curve does \textit{not} represent where the profits are zero. Instead, it represents the marginal profit that the REP can draw from producing hydrogen. Whenever the REP has wind production available, it is guaranteed a minimum profit from hydrogen production. The REP can increase this profit through selling to the electricity market at high prices, or through buying additional electricity for hydrogen production at low prices.  Therefore, an REP with wind production available, will generally draw \textit{additional} profits through market interaction when the electricity price is different from its marginal cost curve. If wind production is not available, it only draws profit from prices lower than its marginal cost.
\end{newtext}

We conclude that to evaluate the profit of a market-bidding REP, it could be essential to model the on the price formation process through, e.g., a DC OPF, as the presence of the electrolyzer in the system can highly impact the LMPs and the REP's profit. This becomes more important with larger electrolyzer capacities and higher hydrogen prices.

\input{images/network_results_table}

\subsection{The Effect of Transmission Network Modeling}
Finally, we investigate the effect of transmission network modeling on the system- and REP-focused results. We compare three network models with varying levels of detail regarding transmission network constraints: a nodal model (which includes all transmission lines), a zonal model (which relaxes all \textit{intra}-regional transmission constraints), and a copper-plate model (which relaxes \textit{all} transmission constraints). To create the zonal model, we define each of the three regions in the RTS-GMLC system, as described by Fig.~\ref{fig:reginal_capacities}, as individual zones. 
The total \textit{inter}-regional transmission capacities between regions 1$\leftrightarrow$2, 1$\leftrightarrow$3, and 2$\leftrightarrow$3, are 1175, 500, and 500 MW, respectively. 
We solve the DC OPF problem for the system with (\textit{i}) no REP (base case) and (\textit{ii}) a market-bidding, 1000 MW electrolyzer REP with \textit{low} hydrogen price, and (\textit{iii}) a market-bidding, 1000 MW electrolyzer REP with \textit{high} hydrogen price. These three cases are solved for each transmission network model, resulting in a total of nine cases.  
For each case, we compute the system load, 
cost of generation per load, system emissions per load, and total renewable power curtailment.
The system-focused results are reported in Table~\ref{tab:transmission_sensitivity}.
We further report REP-focused results including the average electrolyzer load, REP's profit, and the average LMP at the REP's node in Table~\ref{tab:rep_transmission_sens}. All results are computed in the same manner as described in Section~\ref{sec:sys_results}.

\subsubsection{System-focused results} We first examine the system load (first column of Table~\ref{tab:transmission_sensitivity}). For both hydrogen prices, the increase in system load from the base system is lower under a nodal model than under a zonal or copper-plate model. In the case of the higher hydrogen price, the consumption of the REP increases the system load by 23\% under a zonal model and by 17\% under a nodal model. This is expected, as our analysis in Section~\ref{sec:REP_results} showed that the 1000 MW electrolyzer utilization is limited by system congestion and increased LMPs under the nodal model. 
However, there is no substantial further increase in system load when going from a zonal model to a copper-plate model at any hydrogen price. This suggests that congestion is primarily \textit{intra}-regional and can only be fully captured by a nodal model.

Next, we discuss the cost of generation per load (second column). We observe that this cost remains relatively unchanged across the different transmission grid models, even though the nodal model serves less load.  
This indicates that the average cost of generation is higher for the nodal model. This is expected, as congestion in the network, which cannot be captured by a copper-plate model, may lead to the dispatch of more expensive generation in the nodal model. 

The system emissions per load (third column) increase when the transmission network modeling level is relaxed. As seen before, large increases in system load is mainly met by conventional generation, which decreases the share of load met by renewable generation in the system. The system emissions per load are therefore higher for the relaxed transmission network models as they allow for higher REP consumption.

The last column of Table~\ref{tab:transmission_sensitivity} shows the renewable power generation curtailment. These results indicate that the copper-plate and zonal models have much lower curtailment values compared to the nodal model. This is in part because network congestion limits REP consumption and in part because relaxed grid models overestimate the ability to dispatch renewable generation. The copper-plate and zonal network models significantly underestimate the amount of renewable curtailment for the base system and fail to account for the fact that the REP, with a high hydrogen price, increases congestion and subsequent renewable power curtailment.

We conclude that the choice of transmission network model significantly affects the perceived results of adding an electrolyzer to a power system, as the added load could induce congestion in the system. \begin{newtext}
Further, we highlight the role of nodal pricing in mitigating network congestion, as system congestion is captured in this price formation process, as well as the broader benefits of flexible electrolyzer operation. Our results indicate that when the electrolyzer operates flexibly and the REP actively participates in the electricity market, system congestion is either slightly reduced or does not increase significantly. This is reflected in the relatively stable total system cost of generation per unit of load for the market-bidding REPs, as shown in the rightmost plot of Fig.~\ref{fig:sys_costs}. 
\end{newtext}

\input{images/REP_netowrk_results}

\subsubsection{REP-focused results} The results for electrolyzer utilization (first column of Table~\ref{tab:rep_transmission_sens}) further illustrate to what extent the relaxed grid models overestimate the power consumption of the REP. For example, in the case with a higher hydrogen price, the electrolyzer utilization is at 100\% under the relaxed network models, compared to 73\% under the nodal model. This increase in utilization contributes to an overestimation of the REP's profits, which are reported in the second column. 
Further, the increase in perceived profit is in part due to the LMPs at the REP's node (last column of Table~\ref{tab:transmission_sensitivity}). The copper-plate and zonal models yield lower average LMPs, as reduced congestion allows for the dispatch of cheaper generators. In the case of the higher hydrogen price, this leads to vastly differing results on REP profitability. Specifically, the zonal and copper-plate models result in more than twice the REP's profit compared to the nodal model. The results presented for a nodal model in Section~\ref{sec:REP_results} show that a 1000 MW electrolyzer would actually yield a lower profit than a 500 MW electrolyzer in the case of a high hydrogen price, due to the increase in LMPs. This could dis-incentivize larger electrolyzer capacity investments under a nodal pricing scheme. However, the results of profits in Table~\ref{tab:rep_transmission_sens} show that this incentive is lost under the zonal and copper-plate models.

We conclude that relaxed transmission models might result in misleading estimates of the REP's profit as they cannot capture the full effects of congestion on LMPs. Further, these effects lead to investment signals that would be lost under zonal pricing schemes.

%% file: images/REP_results_table.tex
\begin{table}[!t]
\begin{center}
\caption{\small REP results over the year.}
\label{tab:REP_results}
\resizebox{\linewidth}{!}
{
\begin{tabular}{rllrrrrrrrrr} \toprule
    \multicolumn{3}{c}{\textbf{ }}
   &
   \textbf{Electrolyzer}&
  \textbf{REP} &
  \textbf{Average}
  \\ 
  &
   &
   &
  \textbf{average load} &
  \textbf{profit} &
  \textbf{LMP}$^{\rm a}$\\
 \multicolumn{3}{c}{\textbf{Case}} &
  [MW] &
  [M \$] &
  [\$/MWh]\\
  \toprule
\multirow{4}{*}{100 MW}  & \multirow{2}{*}{$\$1.5$/kg} & bidder & 91.44   & 28.88  & 17.81 \\
                         &                             & fixed    & 91.44   & 28.14  & 17.75\\
                         & \multirow{2}{*}{$\$6$/kg}   & bidder & 100.00   & 98.12  & 17.83\\
                         &                             & fixed    & 100.00   & 98.12  & 17.83\\ \midrule
\multirow{4}{*}{500 MW}  & \multirow{2}{*}{$\$1.5$/kg} & bidder & 421.15  & 59.73  & 21.55 \\
                         &                             & fixed    & 421.15  & 47.36  & 24.30\\
                         & \multirow{2}{*}{$\$6$/kg}   & bidder & 499.12 &  362.15 & 31.77\\
                         &                             & fixed    & 499.01 &  352.36 & 33.92\\ \midrule
\multirow{4}{*}{1000 MW} & \multirow{2}{*}{$\$1.5$/kg} & bidder & 629.77 & 67.76  & 24.45\\
                         &                             & fixed    & 617.38 & -730.61 & 179.60\\
                         & \multirow{2}{*}{$\$6$/kg}   & bidder & 729.34 & 320.83 &  81.54\\
                         &                             & fixed    & 672.3 & -563.48  & 241.06 \\ \bottomrule
\end{tabular}}
\end{center}
\footnotesize{$^{\rm a}$ Average LMP at the REP's node.}
\end{table}

%% file: images/network_results_table.tex
\begin{table}[!t]
\caption{\small System-focused results under nodal, zonal, and copper-plate modeling. The results are reported for the base system and two market bidding REPs, both with 1000 MW electrolyzer installed.}
\label{tab:transmission_sensitivity}
{
\centering
\begin{tabular}{llrrrr} \toprule
\textbf{} &
  \textbf{} &
  \multicolumn{1}{c}{\textbf{System}} &
  \multicolumn{1}{c}{\textbf{Cost of}} &
  \multicolumn{1}{c}{\textbf{System}} &
  \multicolumn{1}{c}{\textbf{Renewable}} 
  \\
\textbf{} &
  \textbf{} &
  \multicolumn{1}{c}{\textbf{load}} &
  \multicolumn{1}{c}{\textbf{generation}} &
  \multicolumn{1}{c}{\textbf{emissions }} &
  \multicolumn{1}{c}{\textbf{curtailment}} 
  \\
\multicolumn{2}{c}{\textbf{Case}}  &
  [TWh] &
  [$\$$/MWh] &
  [t/MWh] &
  [GWh] 
  \\ \toprule
\multirow{3}{*}{Base}    & Copper & 37.66    & 18.54      & 0.36   & 348.75     \\
                         & Zonal  & 37.66     & 18.55     & 0.36    & 369.97   \\
                         & Nodal  & 37.66   & 18.67   & 0.37   & 608.09    \\
                         \midrule
\multirow{3}{*}{$\$$1.5/kg} & Copper & 45.07  & 18.81 & 0.43 &  40.47    \\
                         & Zonal  & 44.94  & 18.80 & 0.42 & 40.77     \\
                          & Nodal  & 43.19  & 18.72 &  0.41  & 67.29   \\ 
                         \midrule
\multirow{3}{*}{$\$$6/kg}& Copper & 46.44  & 19.00 & 0.43 & 40.47    \\
                         & Zonal  & 46.44  & 19.04 & 0.43 & 40.77    \\
                         & Nodal  & 44.06  & 19.25 & 0.42  & 302.07  \\ \bottomrule 
\end{tabular}
}
\end{table}

%% file: images/REP_netowrk_results.tex

\begin{table}[t!]
\caption{\small REP-focused results under nodal, zonal, and copper-plate modeling. Results are reported for a 1000 MW electrolyzer, market bidding REP, under the lower and higher hydrogen price. }
\label{tab:rep_transmission_sens}
\centering
\begin{tabular}{llrrr} \toprule
&                 
& \multicolumn{1}{c}{\textbf{Electrolyzer}} 
& \multicolumn{1}{c}{\textbf{REP}} 
& \multicolumn{1}{c}{\textbf{Average}} \\
&                 
& \multicolumn{1}{c}{\textbf{average load}} 
& \multicolumn{1}{c}{\textbf{profit}}
& \multicolumn{1}{c}{\textbf{LMP}}\\
\multicolumn{2}{c}{\textbf{Case}} 
& { [}MW{]}   
& {[}M $\$${]}                 
& { [}\$/MWh{]}           \\ \toprule
\multirow{3}{*}{$\$$1.5/kg} & Copper & 843.53  & 80.92   & 22.05 \\ 
                            & Zonal  & 829.25  & 79.80   & 22.26 \\
                            & Nodal  & 629.77  & 67.76   & 24.45 \\ \midrule
\multirow{3}{*}{$\$$6/kg}   & Copper & 1000.00    & 769.76  & 22.39 \\
                            & Zonal  & 1000.00    & 766.06  & 22.85 \\
                            & Nodal  & 729.34  & 320.83  & 81.55 \\ \bottomrule
\end{tabular}
\end{table}

%% file: 4_Conclusion.tex
\section{Conclusions}
This paper introduces the derivation of the opportunity-cost curve of REPs, with the aim of facilitating the integration of such plants in electricity market-clearing models and the analysis of their system-wide impacts on emissions, costs, and congestion, under perfect competition. The opportunity cost is derived from the revenue of hydrogen production and forms a curve that can be used to model REP bids in an electricity market. We compare a market-bidding REP, bidding its derived opportunity-cost curve, to a REP with a fixed electrolyzer consumption. The results are obtained by solving a DC OPF of the RTS-GMLC test system, modifying a pre-existing wind farm cost curve to represent the REP. 

The results show that active market participation can efficiently mitigate increases in the cost of generation per load that occur when the electrolyzer is operated in a fixed manner. \begin{revised} A market-bidding REP also efficiently reduces the curtailed renewable power in the system compared to the fixed electrolyzer case. However, high hydrogen prices will move the REP further up on the merit-order curve. This means that the REP's electrolyzer will consume power at higher cleared electricity prices, which causes it to behave more like a fixed load in the system. Therefore, higher hydrogen prices reduces the system-level benefits of a market-bidding REP. Further,\end{revised} market-bidding in itself is not necessarily an efficient measure to prevent increases in system emissions due to the lack of correlation between the emission intensities and the marginal cost of generators. 

The size of the electrolyzer relative to the local renewable capacity and surrounding transmission network is of great importance, as a large electrolyzer capacity might induce congestion in the system, amplifying the negative effects on the cost of generation and the renewable power curtailment. Modeling the full transmission network when evaluating the presence of an electrolyzer in the power system might therefore be important because the effects of a congested network are otherwise unknown. 
However, if the REP is price sensitive both in theory (bidding its marginal cost curve) and in practice (marginal cost proportional to the system), these negative effects are largely mitigated.  

This work leaves several future research directions. The REP might be exposed to uncertainty in its renewable production as well as to hydrogen demand constraints, which will affect the opportunity cost of selling electricity. Further, the REP's potential as a price maker, and to what extent it might act strategically rather than bidding its true opportunity cost curve, can be analyzed. 
\begin{revised}
    The developed REP bid curve could be used in an analysis to evaluate the long-term impacts on, e.g., system emissions.
\end{revised} 
Finally, this work can be extended to investigate the impact of renewable hydrogen regulations and subsidies on system and REP results.

%% file: 5_appendix.tex
\subsection{Derivation of the Approximated Cost Curve Slopes and Intercepts} 
\label{app:cost_curve}
\begin{newtext}
For notational simplicity, we denote $x$ as the general  variable for power in all functions through the appendix. The empirical hydrogen production function is first approximated as a piecewise linear function, as given by eq.~(\ref{eq:hydrogen_piecewise}):
\begin{align}
\tag{\ref{eq:hydrogen_piecewise}}
 h(x) = \sum_{i \in \mathcal{I}} \left(A_i x + B_i\right)\mathbf{1}_{\left(x \in [\underline{P}_i, \overline{P}_i ]\right)},
\end{align}
where $i$ is the piece of the approximation, \( \mathcal{I} \in \{1,2, \ldots, I \} \), $A_i$ and $B_i$ are the corresponding slopes and intercepts respectively, and $\underline{P}_i, \overline{P}_i$ are the breakpoints:
\[\underline{P}_i < \overline{P}_{i} = \underline{P}_{i+1}<... <\overline{P}_{I},\]  
where $\overline{P}_{I} = \bar{P}^{\rm h}$ and $\underline{P}_{0} = 0$.
We define the approximated  hydrogen \textit{revenue} function based on eq.  (\ref{eq:hydrogen_revenues}) and eq.~(\ref{eq:hydrogen_piecewise}) as

\begin{align}
\label{eq:rev_piecewise}
{\rm r}(x) = \lambda^{\rm{h}}
   \begin{cases}
      A_i x + B_i, \quad x \in [\underline{P}_i, \overline{P}_i],  \quad i \in \mathcal{I}\\
        A_{i} \overline{P}_i + B_{i}, \quad x > \overline{P}_i, \quad \quad \quad    i = I.
   \end{cases}
\end{align}

Hence, if $x$ exceeds the electrolyzer capacity, i.e., if the wind power production exceeds the electrolyzer capacity, the piecewise revenue function (\ref{eq:rev_piecewise}) will have one more piece than the piecewise hydrogen function, $I + 1$.
In order to simplify eq.~(\ref{eq:rev_piecewise}), we represent the final piece on the revenue curve through new elements on the $\bm A$ and $\bm B$ vectors. We expand the set $\mathcal{I}$, so that $\mathcal{I} = \{1,\ldots,I,I+1\}$. The additional piece on eq.~(\ref{eq:rev_piecewise}) is simply a constant. We achieve $A_{I+1} = 0$ and $B_{I+1} = A_{I} \overline{P}_I + B_{I}$. The final breakpoint of the curve, now defined by $\overline{P}_{I+1}$, is for simplicity set to $\overline{P}_{I+1} = \bar{P}^{\rm RES}$, where $\bar{P}^{\rm RES}$ is the maximum RES capacity. The approximated hydrogen revenue eq.~(\ref{eq:rev_piecewise}) is then expressed as:
\begin{align}
{\rm r}(x) = \lambda^{\rm{h}}(
      A_i x + B_i), && x \in [\underline{P}_i, \overline{P}_i],  && i \in \mathcal{I}.
\end{align}

Recall that the REP's cost function, as expressed in eq.~(\ref{eq:op_cost_curve}), is equal to the difference between the revenue function evaluated at two different points,
\begin{align}
    {\rm c}(x) =   r(P^{{\rm RES}}) - r(P^{\rm RES} - x) \tag{\ref{eq:op_cost_curve}},
\end{align}
where \( P^{\rm RES} \) is a fixed scalar, i.e., the first term is a constant. 

We identify the piece corresponding to \( P^{\rm RES} \), finding the index \( j \) such that
\begin{align}
    P^{\rm RES} \in [\underline{P}_i, \overline{P}_i].
\end{align}
We can then express the first term of the REP's cost function through the slope $A$ and intercept $B$ of piece $i$ of the approximated hydrogen production curve,
\begin{align}
    {\rm r}(P^{\rm RES}) = \lambda^{\rm h}(A_j P^{\rm RES} + B_j), && P^{\rm RES} \in [\underline{P}_j, \overline{P}_j], && j \in \mathcal{I}.
\end{align}

Next, we evaluate the second term of the REP's cost curve eq.~(\ref{eq:op_cost_curve}). This term is piecewise linear in the variable \( x \),
\begin{align}
{\rm r}(P^{\rm RES} - x) = \lambda^{\rm h} ( A_j (P^{\rm RES} - x) + B_j), 
\end{align}
for, 
\begin{align}
    \label{eq:j_piece_selection}
    (P^{\rm RES} - x) \in [\underline{P}_i, \overline{P}_i].
\end{align}
The piece range of eq.~(\ref{eq:j_piece_selection}) can be rewritten as,
\[
x \in [P^{\rm RES} - \overline{P}_i, \; P^{\rm RES} - \underline{P}_i],
\]
which allows us to express the second term of the REP's cost curve based on the approximation of the hydrogen production curve. We get,
\begin{align*}
    r(P^{\rm RES} - x) 
    &= \lambda^{\rm{h}} ( A_i P^{\rm RES} + B_i - A_i x ),\\
    & x \in [P^{\rm RES} - \overline{P}_i , P^{\rm RES} - \underline{P}_i].
\end{align*}

With both terms of the REP cost curve eq.~(\ref{eq:op_cost_curve}) expressed as piecewise linear approximations, we derive: 
\begin{align}
    {\rm c}(x) &= {\rm r}(P^{\rm RES}) - {\rm r}(P^{\rm RES} - x) \notag \\
&= \lambda^{\rm{h}} \bigl(A_j P^{\rm RES} + B_j \bigr) - \lambda^{\rm{h}} \bigl(A_i P^{\rm RES} + B_i - A_i x \bigr) \notag \\
&= \underbrace{\lambda^{\rm{h}}\bigl( (A_j - A_i) P^{\rm RES} + (B_j - B_i) \bigr)}_{\text{intercept}} + \underbrace{\lambda^{\rm{h}}A_i}_{\text{slope}}  x,
\end{align}
which give us the expressions of the slope $\alpha_i$ and intercept $\beta_i$ of the approximated cost curve.
\begin{align}
{\rm c}(x) 
&= \alpha_ix+\beta_i, \\
\alpha_i &= \lambda^{\rm{h}}A_i, \\
\beta_i &= \lambda^{\rm{h}}\bigl( (A_j - A_i) P^{\rm RES} + (B_j - B_i) \bigr), 
\end{align}
where $x  \in [P^{\rm RES} - \overline{P}_i , P^{\rm RES} - \underline{P}_i], \quad i \in \mathcal{I} = \{0,1,\ldots, I+1\}$ and $j$ is the index found so that $P^{\rm RES} \in [\underline{P}_j,\overline{P}_j], \quad j \in \mathcal{I} $.
\end{newtext}